\documentclass[12pt]{amsart}
\usepackage{amscd,amsfonts,amsmath,amsthm,amssymb,verbatim,mathrsfs}
\usepackage[dvips]{graphicx}
\usepackage{tikz}
\usetikzlibrary{shapes,backgrounds}
\usetikzlibrary{calc}
\usetikzlibrary{matrix}
\usepackage{graphics,hyperref}
\usepackage{pstricks}
\usepackage{enumitem}
\usepackage{psfrag}
\usepackage{soul}
\usepackage{pstcol,pst-plot,pst-3d}
\usepackage{latexsym}
\usepackage{pstcol,pst-plot,pst-3d}
\usepackage{multicol}
\usepackage[normalem]{ulem} 

\psset{unit=0.7cm,linewidth=0.8pt,arrowsize=2.5pt 4}

\newpsstyle{fatline}{linewidth=1.5pt}
\newpsstyle{fyp}{fillstyle=solid,fillcolor=verylight}
\definecolor{verylight}{gray}{0.97}
\definecolor{light}{gray}{0.9}
\definecolor{medium}{gray}{0.85}

\newtheorem{thm1}{Theorem}[section]

\newtheorem{rem1}[thm1]{Remark}

\newtheorem{prop1}[thm1]{Proposition}

\begin{document}

\title{On the relative size of toric bases}
\author{Christos Tatakis and Apostolos Thoma}
\address{Christos Tatakis, Department of Mathematics, University of Ioannina, Ioannina 45110, Greece}
\email{chtataki@uoi.gr}

\address{Apostolos Thoma, Department of Mathematics, University of Ioannina, Ioannina 45110, Greece} 
\email{athoma@uoi.gr}

\begin{abstract}
We consider the Graver basis, the universal Gr\"{o}bner basis, a Markov basis and the set of the circuits of a toric ideal. 
Let $A, B$ be any two of these bases such that $A\not \supseteq B$ in general, we prove that there is no polynomial on the size  or
on the maximal degree of the elements of $A$  which bounds the size or the maximal degree of the elements of $B$ correspondingly. 

\end{abstract}
\maketitle

\section{Introduction}

\par  Let $A=\{\textbf{a}_1,\ldots,\textbf{a}_m\}\subseteq \mathbb{N}^n$
be a nonzero vector configuration in $\mathbb{Q}^n$ and
$\mathbb{N}A:=\{l_1\textbf{a}_1+\cdots+l_m\textbf{a}_m \ | \ l_i \in\mathbb{N}\}$ the 
corresponding affine semigroup. There are two cases for the semigroup  $\mathbb{N}A$
either it
is pointed, that is $\mathbb{N}A\cap (-\mathbb{N}A)=\{0\}$, or it is not pointed.
  We grade the
polynomial ring $\mathbb{K}[x_1,\ldots,x_m]$ over an arbitrary field $\mathbb{K}$ by the
semigroup $\mathbb{N}A$ setting $\deg_{A}(x_i)=\textbf{a}_i$ for
$i=1,\ldots,m$. For $\textbf{u}=(u_1,\ldots,u_m) \in \mathbb{N}^m$,
we define the $A$-\emph{degree} of the monomial $\textbf{x}^{\textbf{u}}:=x_1^{u_1} \cdots x_m^{u_m}$
to be \[
\deg_{A}(\textbf{x}^{\textbf{u}}):=u_1\textbf{a}_1+\cdots+u_m\textbf{a}_m
\in \mathbb{N}A,\] 
while we denote the usual degree $u_1+\cdots +u_m$ of 
$\textbf{x}^{\textbf{u}}$ by $\deg(\textbf{x}^{\textbf{u}})$. The
\emph{ toric ideal} $I_{A}$ associated to $A$ is the prime ideal
generated by all the binomials $\textbf{x}^{\textbf{u}}- \textbf{x}^{\textbf{v}}$
such that $\deg_{A}(\textbf{x}^{\textbf{u}})=\deg_{A}(\textbf{x}^{\textbf{v}})$, see \cite{St}.\par 
There are several sets for a toric ideal which include
crucial information about it, such as the Graver basis,
the Markov bases, the universal Gr\"obner basis and
the set of the circuits.
An irreducible binomial $\textbf{x}^{\textbf{u}}-
\textbf{x}^{\textbf{v}}$ in $I_A$ is called primitive if
there is no other binomial
 $\textbf{x}^{\textbf{w}}- \textbf{x}^{\textbf{z}}$ in $I_A$,
such that $\textbf{x}^{\textbf{w}}$ divides $
\textbf{x}^{\textbf{u}}$ and $\textbf{x}^{\textbf{z}}$ divides $
\textbf{x}^{\textbf{v}}$. The set of the primitive binomials is finite, forms the
Graver basis of $I_A$ and is denoted by $Gr_A$.   
 The universal Gr\"{o}bner basis 
 of an ideal $I_A$ is defined as the union of all reduced Gr\"obner bases $G_<$ of $I_A$, as $<$ runs over all term orders. It is a finite subset of the $I_A$ and it is a Gr\"obner basis for the ideal with respect to all admisible term orders, see \cite{St}. The support of
 a monomial ${\bf{x^u}}$ of $\mathbb{K}[x_1,\ldots,x_m]$
is $\textrm{supp}({\bf{x^u}}):=\{i\ | \ x_{i}\ \mbox{ divides}\ {\bf{x^u}}\}$ and
the support of a binomial $B={\bf{x^u}}-{\bf{x^v}}$ is
$\textrm{supp}(B):=\textrm{supp}({\bf{x^u}})\cup \textrm{supp}({\bf{x^v}})$.
An irreducible nonzero binomial is called circuit if it has minimal support.
The set of the circuits of a toric ideal $I_A$ is denoted by
$\mathcal{ C}_A$. A  Markov basis is a minimal generating set of the toric ideal $I_A$, 
consisting of binomials, see \cite[Theorem 3.1]{DST}. A relation between some of the above sets was given by B. Sturmfels
in \cite{St}:

\begin{prop1}\label{CUG}\cite[Proposition 4.11]{St} For any toric ideal $I_A$ it holds:
$$\mathcal{ C}_A\subseteq \mathcal{ U}_A \subseteq  Gr_A.$$
\end{prop1} 

Every reduced Gr\"{o}bner basis is a generating set of the toric ideal $I_A$ consisting of binomials, therefore
it contains also a Markov basis. Thus the universal Gr\"{o}bner basis and the Graver basis contain at least one Markov basis.
  The Graver basis contains all the
 Markov bases of $I_A$ if and only if the semigroup $\mathbb{N}A$ is pointed, see \cite{THO2}. 
 It is well known that the above inclusions may or may not be strict, see \cite[Example 4.12]{St}. In famous classes of ideals, the equality
 happens between some of the above bases and is combined with interesting geometric, combinatorial and homological properties.
 For example robust are those  toric ideals for which the 
 universal Gr\"{o}bner basis is a Markov basis \cite{BR} and strongly robust are those toric ideals 
 for which the Graver basis is a Markov basis, see \cite{PeThVl, Sull}. Lawrence toric ideals are strongly robust \cite{St} but also 
 toric ideals of non pyramidal self dual projective toric varieties 
 are strongly robust \cite{TV}. For unimodular toric ideals all the elements in the Graver basis are circuits \cite{BPS}.

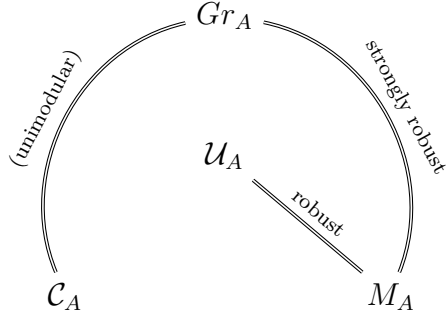
\begin{figure}[h]
\begin{center}
\begin{tikzpicture}[descr/.style={fill=white,inner sep=1.5pt}]
    \matrix (m) [matrix of math nodes,row sep=3em,column sep=3em]
  {
       & Gr_A &  \\
      & \mathcal{U}_A &  \\  
			\mathcal{C}_A &  & M_A	\\		
		};
     \path[-stealth]
    (m-1-2) (m-2-2)
 (m-2-2) edge [double,-] node[sloped, above] {\tiny{robust}}(m-3-3) 
						 (m-1-2)  edge [double,-,bend left=50] node[sloped, above] {\tiny{strongly robust}} (m-3-3)
						 (m-3-1)edge [double,-,bend left=50] node[sloped, above] {\tiny{(unimodular)}} (m-1-2)
						 (m-3-1)  (m-3-3)
						 (m-2-2)  (m-3-1);
\end{tikzpicture}
\caption{Well-known classes of toric ideals}
\label{Fig.well-known-classes}
\end{center}
\end{figure}

There are several results in the literature concerning degree bounds of the elements of these sets and sometimes bounds 
on the one of these sets in terms of another set. There exist several bounds 
on the degrees of the elements of the Graver basis of a toric ideal which 
have important implications to integer programming and computational algebraic geometry,
see for example \cite{LST, H, Pe, RYCH, StTat, St, St1, TT1}.

The aim of this article is to present several theorems concerning bounds on the size of these bases or the maximal degree of their elements  in terms 
of the size or the maximal degree of the other bases. The proofs are based on carefully chosen simple counterexamples to show these
relations with the fewest possible examples. Some of these examples are well known and we include them for completeness. There are also several other examples 
showing this extremal behavior. 

In Section 2 we present the basic results 
about the toric ideals of graphs which will be useful for us in the sequel. For more details 
we refer to \cite{RTT, TT1, Vi}.

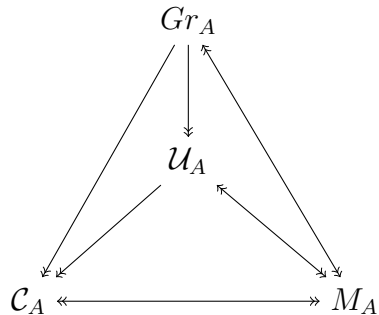
\begin{figure}[h]
\begin{center}
\begin{tikzpicture}[descr/.style={fill=white,inner sep=1.5pt}]
    \matrix (m) [matrix of math nodes,row sep=3em,column sep=3em]
  {
       & {Gr_A}  &  \\
      & {\mathcal{U}_A} &  \\  
			{\mathcal{C}_A}&  & {M_A}  \\		
		};
     \path[-stealth]
    (m-1-2) edge [->>] (m-2-2)
						 (m-3-3) edge[<<->>] (m-2-2) 
						 (m-1-2) edge[<<->>] (m-3-3)
						 (m-3-1) edge[<<-]  (m-1-2)
						 (m-3-1) edge[<<->>] (m-3-3)
						 (m-2-2) edge[->>]  (m-3-1);
\end{tikzpicture}
\caption{The relations of the size and the maximal degrees of the elements of the toric bases}
\label{relations!}
\end{center}
\end{figure}

In Section 3 are the main results of the article that can be summarized in Figure \ref{relations!}. 
 In Figure \ref{relations!}, $B\twoheadrightarrow C$ represents that the size of the base $B$ or the degrees of the elements of the set $B$ cannot be bounded above by a polynomial on the size or the maximal degree of the elements of $C$.

\section{On the toric bases of toric ideals of graphs}\label{section 2}

\par

Let $G$ be a finite simple connected graph with vertices
$V(G)=\{v_{1},\ldots,v_{n}\}$ and edges $E(G)=\{e_{1},\ldots,e_{m}\}$.
Let $\mathbb{K}[e_{1},\ldots,e_{m}]$
be the polynomial ring in the $m$ variables $e_{1},\ldots,e_{m}$ over a field $\mathbb{K}$.  We
will associate each edge $e=\{v_{i},v_{j}\}\in E(G)$ with the element
$a_{e}=v_{i}+v_{j}$ in the free abelian group $ \mathbb{Z}^n $
with basis the set of the vertices
of $G$.  Each vertex $v_j\in V(G)$ is associated with the vector $(0,\ldots,0,1,0,\ldots,0)$, where the 
nonzero component is in the $j$ position. We denote by $I_G$ the toric ideal $I_{A_{G}}$ in
$\mathbb{K}[e_{1},\ldots,e_{m}]$, where  $A_{G}=\{a_{e}\ | \ e\in E(G)\}\subset \mathbb{Z}^n $.

A \emph{walk}  connecting $v_{i_{1}}\in V(G)$ and
$v_{i_{s+1}}\in V(G)$ is a finite sequence of the form
$$w=(\{v_{i_1},v_{i_2}\},\{v_{i_2},v_{i_3}\},\ldots,\{v_{i_s},v_{i_{s+1}}\})$$
with each $e_{i_j}=\{v_{i_j},v_{i_{j+1}}\}\in E(G)$, for $j=1,\ldots,s$. A \emph{trail} (respectively \emph{path}) is a walk in which all edges (respectively vertices) are distinct.
The \emph{length}
of the walk $w$ is  the number $s$ of its edges. An
even (respectively odd) walk is a walk of \emph{even} (respectively odd) length.
A walk
$w=(\{v_{i_1},v_{i_2}\},\{v_{i_2},v_{i_3}\},\ldots,\{v_{i_s},v_{i_{s+1}}\})$
is called \emph{closed} if $v_{i_{s+1}}=v_{i_1}$. A \emph{cycle}
is a closed walk
$(\{v_{i_1},v_{i_2}\},\{v_{i_2},v_{i_3}\},\ldots,\{v_{i_s},v_{i_{1}}\})$ with
$v_{i_k}\neq v_{i_j},$ for every $ 1\leq k < j \leq s$.

Given an even closed walk $w$ of the graph $G$, where $w =(e_{i_1}, e_{i_2},\dots,
e_{i_{2q}}),$ we denote by $B_w$ the binomial
$$B_w=E^+(w)-E^-(w),$$
where $E^+(w)=\prod _{k=1}^{q} e_{i_{2k-1}},\ E^-(w)=\prod _{k=1}^{q} e_{i_{2k}}$.
It is known that the toric ideal $I_G$
is generated by binomials of this form, see \cite{Vi}. Note that the binomials $B_w$ 
are homogeneous and the degree of $B_w$ is $q$, the half of the number of the edges of the walk. 
For convenience,
we denote by $\textbf{w}$ the subgraph of $G$ with vertices the vertices of the
walk and edges the edges of the walk $w$. We call a walk
$w'=(e_{j_{1}},\dots,e_{j_{t}})$ a \emph{subwalk} of $w$ if
$e_{j_1}\cdots e_{j_t}| e_{i_1}\cdots e_{i_{2q}}.$
An even
closed walk $w$ is said to
be primitive if there exists no even closed subwalk $\xi$ of $w$ of smaller
length such that $E^+(\xi)| E^+(w)$ and $E^-(\xi)| E^-(w)$. The walk $w$
is primitive if and only if the binomial $B_w$ is primitive.

A \emph{ cut edge} (respectively \emph{ cut vertex}) is an edge (respectively vertex) of
the graph whose removal increases the number of connected
components of the remaining subgraph.  A graph is called \emph{
biconnected} if it is connected and does not contain a cut
vertex. A \emph{ block} is a maximal biconnected subgraph of a given
graph $G$.

The following theorems determine the form of the circuits and the
primitive binomials of a toric ideal of a graph $G$. R.~Villarreal in
\cite[Proposition 4.2]{Vi} gave a necessary and sufficient
characterization of the circuits:
\begin{thm1}\label{circuit} Let $G$ be a graph and let $W$ be a connected subgraph of $G$.
The subgraph $W$ is the graph  ${\bf w}$ of a walk $w$ such that  $B_w$ is a circuit
 if and only if
\begin{enumerate}
  \item $W$ is an even cycle or
  \item $W$ consists of two odd cycles intersecting in exactly one vertex or
  \item $W$ consists of two vertex-disjoint odd cycles joined by a path.
\end{enumerate}
\end{thm1}
 
The next theorem by  E.~Reyes et all, see \cite{RTT}, describes  the form of the underlying graph
of a primitive walk and thus describes the Graver basis of $I_G$.

\begin{thm1} \label{primitive-graph}
Let $G$ be a  graph and let $W$ be a connected subgraph of $G$.
The subgraph $W$ is the graph  ${\bf w}$ of a primitive walk $w$
 if and only if \begin{enumerate}
  \item  $W$ is an even cycle or
  \item  $W$ is not biconnected and
\begin{enumerate}
  \item every block of $W$ is a cycle or a cut edge and
  \item every cut vertex of $W$ belongs to exactly two blocks and separates the graph in two parts, 
  the total number of edges
of the  blocks that are cycles in each part is odd.
\end{enumerate}
\end{enumerate}
\end{thm1}

We remark that every even primitive walk $w=(e_{i_1},\ldots,e_{i_{2k}})$
partitions the set of the edges in the two sets $w^+= \{e_{i_j}|j \
{\it odd}\}$ and $w^-=\{e_{i_j}|j \ {\it even}\}$, otherwise the
binomial $B_w$ is not irreducible. The edges of $w^+$ are called odd edges of the walk w and those of
$w^-$ are called even. If $e_i\in w^+$ and $e_j\in w^-$, we say that the edges $e_i$ and $e_j$ of the walk $w$ have different parity. \emph{Sink} of a block $B$ is a
common vertex of two odd or two even edges of the walk $w$ which
belong to the block $B$.
The last condition of Theorem \ref{primitive-graph} 
can be expressed also in terms of the walk as: 
every cut vertex of $\bf{w}$ belongs to exactly two blocks and it is a sink of both. A sink of a block
should be always a cut vertex \cite{RTT}. 

Afterwards, we recall from \cite{RTT}, some graph theoretical notions
 in order to describe a Markov basis of a toric ideal of a graph $G$.
 A binomial is called minimal,
if it belongs to at least one minimal system of generators of $I_G$, i.e.
at least one Markov basis of $I_G$.

For a given subgraph $F$  of $G$, an edge $f$
of the graph $G$ is called chord of the subgraph $F$, if the vertices of the edge
$f$ belong to $V(F)$ and $f\notin E(F)$. A chord $e=\{v_k,v_l\}$ is called
bridge of a primitive walk $w$ if there exist two different blocks
${\mathcal B}_1,{\mathcal B}_2$ of $\bf{w}$ such that $v_k\in {\mathcal B}_1$
and $v_l\in {\mathcal B}_2$. Let $w$ be an even closed walk
$(\{v_{1},v_{2}\},\{v_{2},v_{3}\},\ldots,\{v_{2q},v_{1}\})$ and
$f=\{v_{i},v_{j}\}$ a chord of $w$. Then, $f$ breaks $w$ into two
walks:
$$w_{1}=(e_{1},\ldots,e_{i-1}, f, e_{j},\ldots,e_{2q})\ \textrm{and}\
w_{2}=(e_{i},\ldots,e_{j-1},f),$$ where $e_{s}=\{v_{s},v_{s+1}\},\ 1\leq s< 2q$ and $e_{2q}=\{v_{2q},v_{1}\}.$
The two walks are either both even or both odd. A chord is called even (respectively odd) if it is not
a bridge and it breaks the walk into two even walks (respectively odd). A primitive walk $w$ is called strongly primitive if it has
not two sinks with distance one in any cyclic block of {\bf{w}}.

The next theorem by Reyes et al,
gives a necessary and sufficient
characterization of the minimal binomials of a toric ideal of a graph $G$, thus describes the elements of 
Markov bases of $I_G$.

\begin{thm1}\cite[Theorem 4.13]{RTT} \label{minimal}
Let $w$ be an even closed walk. $B_{w}$ is a minimal binomial if
and only if \begin{enumerate}
  \item[(M1)] all the chords of $w$ are odd,
  \item[(M2)] there are not two odd chords of $w$ which cross effectively except if they form an $F_4$,
  \item[(M3)] no odd chord crosses an $F_4$ of the walk $w$,
  \item[(M4)] $w$ is a strongly primitive.
\end{enumerate}
\end{thm1}

The following theorem determines the indispensable elements of the ideal $I_G$, indispensable are the elements
that belong to every Markov basis of $I_G$.

\begin{thm1}\cite[Theorem 4.14]{RTT}\label{indispensable}
Let $w$ be an even closed walk. $B_{w}$ is an
indispensable binomial if and only if $w$ is a strongly primitive walk, all the chords of $w$ are
odd and there are not two of them which cross effectively.
\end{thm1}

Finally, in order to describe the universal
Gr\"obner basis for the case of toric ideals of graphs, we remind the notions of pure blocks
 and of the mixed walks of a graph $G$, see \cite{TT1}.
 A cyclic block $\mathcal{B}$ of a primitive walk $w$ is called pure if all the edges
of the block $\mathcal{B}$ belong either to $\textbf{w}^+$ or to $\textbf{w}^-$. A primitive
walk $w$ is called mixed if none of the cyclic blocks of $w$ is pure. The next theorem
describes completely the elements of the universal Gr\"obner
basis of a toric ideal of a graph $G$.

\begin{thm1}\cite[Theorem 3.4]{TT1} \label{UGB} Let $w$ be a primitive walk.
$B_w$ belongs to the universal Gr\"{o}bner basis of $I_G$ if and only if $w$ is mixed.
\end{thm1}

\section{Size and degree bounds on toric bases}

\subsection{On the comparison of the size of toric bases}

Proposition \ref{CUG} gives  the comparison of the sets of $\mathcal{C}_A, \mathcal{U}_A, Gr_A $ of a toric ideal. In the case of toric ideals of graphs, we know that every Markov basis of the ideal belongs also to its universal Gr\"obner basis, 
see \cite[Proposition 3.3.]{T}, but this is not true in the general case, see \cite[Example 1.8.]{THO1}. In general, we have no information about the differences of the sizes between the above sets.

It is reasonable to ask about the comparison of the size of a Markov basis of the ideal with the subsets $\mathcal{C}_A, \mathcal{U}_A, Gr_A $ of $I_A$. In this subsection, we answer the above question, studying the problem in two cases, i.e. the case that $\mathbb{N}A$ is pointed and the case that it is not.

\subsubsection{{\bf{Markov basis comparing to the sets of  $\mathcal{C}_A, \mathcal{U}_A, Gr_A $ of an ideal , where $\mathbb{N}A$ is pointed}}}

\begin{figure}[h]
\begin{center}
\psfrag{A}{$v_1$}\psfrag{B}{$v_2$}\psfrag{C}{$v_{n+1}$}\psfrag{D}{$u_1$}
\psfrag{H}{$u_2$}\psfrag{J}{$u_{n+1}$}\psfrag{E}{$\cdots$}\psfrag{M}{$s_1$}\psfrag{N}{$t_1$}\psfrag{O}{$t_n$}
\psfrag{Q}{$x_1$}\psfrag{R}{$y_1$}\psfrag{S}{$y_n$}
\includegraphics{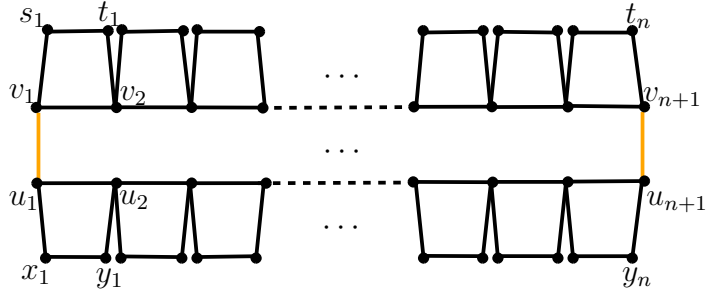}
\caption{Example of a graph $G$ such that $\mid M_G \mid  \ll  \mid \mathcal{C}_G \mid = \mid  \mathcal{U}_G \mid = \mid Gr_G \mid $. }
\label{Markov.vs.others}
\end{center}
\end{figure}

Let $G$ be the graph of the Figure \ref{Markov.vs.others} 
with vertices $V(G)=\{v_i, u_i\ |\ 1\leq i\leq n+1\}\cup \{s_i, t_i, x_i, y_i\ |\ 1\leq i\leq n\}$ and edges 
 \begin{eqnarray}
E(G) &=&\{\{v_i, s_i\}, \{s_i, t_i\}, \{t_i, v_{i+1}\}, \{v_i, v_{i+1}\},
\{u_i, x_i\}, \{x_i, y_i\}, \{y_i, u_{i+1}\}, \{u_i, u_{i+1}\}\ |\ 1\leq i\leq n\} \nonumber \\
&\cup &\{\{v_1, u_1\}, \{v_{n+1}, u_{n+1}\}\}. \nonumber 
\end{eqnarray}

and let $I_G$ be its corresponding toric ideal. 

\begin{figure}[h]
\begin{center}
\begin{tikzpicture}[descr/.style={fill=white,inner sep=1.5pt}]
    \matrix (m) [matrix of math nodes,row sep=3em,column sep=3em]
  {
       & \mid Gr_A\mid &  \\
      & \mid \mathcal{U}_A \mid&  \\  
			\mid \mathcal{C}_A \mid&  & \mid M_A	\mid\\		
		};
     \path[-stealth]
    (m-1-2) (m-2-2)
						 (m-3-3) edge[<<-] (m-2-2) 
						 (m-1-2) edge[->>] (m-3-3)
						 (m-3-1) (m-1-2)
						 (m-3-1) edge[->>] (m-3-3)
						 (m-2-2)  (m-3-1);
\end{tikzpicture}
\caption{Size comparison between a Markov basis and all others}
\label{relations2}
\end{center}
\end{figure}
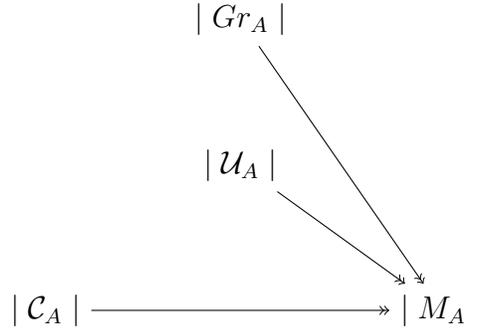

We use the above graph to prove the following theorem, see Figure \ref{relations2}:

\begin{thm1}\label{size1} The size of the elements of the Graver basis, the universal Gr\"obner basis 
and the set of the circuits of a toric ideal $I_A$ cannot
be bounded above by a polynomial on the size of  a Markov basis of $I_A$.
\end{thm1}
\begin{proof} For the graph of the Figure \ref{Markov.vs.others} we claim that $$\mid M_G \mid = 2n+1\ \textrm{and}\  \mid \mathcal{C}_G \mid = \mid  \mathcal{U}_G \mid = \mid Gr_G \mid = 2n+4^n.$$

Note that the graph is bipartite with bipartition $(X,Y)$, where

 $$X=\{s_i,u_i,y_i,t_j,x_j,v_j\ |\ i\ \textrm{odd},\ j\ \textrm{even}\}$$ and $$Y=\{s_j,u_j,y_j,t_i,x_i,v_i\ |\ i\ \textrm{odd},\ j\ \textrm{even}\}.$$ 

It follows that the minimal generators of the corresponding toric ideal are exactly the binomials whose corresponding walks are 
the cycles of the graph $G$ with no chords. There are $2n$ minimal generators of  length four, in the forms $(v_i, s_i, t_i, v_{i+1})$ or 
$(u_i, x_i, y_i, u_{i+1})$  and one 
of length $2n+2$, the cycle $(v_1, v_2, \cdots, v_{n+1}, u_{n+1}, u_n,\cdots, u_2, u_1)$.   It follows that $\mid M_G \mid = 2n+1$.

Toric ideals of bipartite graphs are unimodular thus every element of the Graver basis of the ideal $I_G$ is also a circuit. 
By Proposition \ref{CUG} it follows that $ \mid \mathcal{C}_G \mid = \mid  \mathcal{U}_G \mid = \mid Gr_G \mid $. 
There are $2n$ even cycles of length 4 and  all other cycles are passing through both  edges $\{v_1, u_1\}, \{v_{n+1}, u_{n+1}\}$.
From $v_1$ to $v_{n+1}$ there are $2^n$ different paths and from $u_{n+1}$ to $u_1$ there are also $2^n$ different paths.
Therefore there are $2^n\cdot 2^n=4^n$ cycles passing through the  edges $\{v_1, u_1\}, \{v_{n+1}, u_{n+1}\}$.
Thus $ \mid \mathcal{C}_G \mid = \mid  \mathcal{U}_G \mid = \mid Gr_G \mid = 2n+4^n$. 

In this case the semigroup $\mathbb{N}A_G$ is pointed, thus all Markov bases have the same size, actually in this example there exists only one
Markov basis. Let $s=2n+1$ be the size of the Markov basis then $ \mid \mathcal{C}_G \mid = \mid  \mathcal{U}_G \mid = \mid Gr_G \mid = s-1+4^{(s-1)/2}$. Therefore in this example
the size of the elements of the Graver basis, the universal Gr\"obner basis
and the set of the circuits is exponential on the size of a Markov basis of $I_G$. The result follows. 
\end{proof}

\begin{rem1}{\rm
Note that for this example although the sizes of $\mathcal{C}_G, \mathcal{U}_G, Gr_G$ cannot  be
bounded above by a polynomial on the size of  a Markov basis of $I_A$, their maximal degrees have a linear relation, since
$\deg_{M_A}=n+1$ and $\deg_{Gr_A}=3n+1$. Therefore, in this example the degrees of $\mathcal{C}_G, \mathcal{U}_G, Gr_G$ can be bounded by a linear polynomial on $\deg_{M_A}$. Later on, we prove that this is not true in the general case, see Theorem \ref{TM4}. }
\end{rem1}

\subsubsection{{\bf{Markov basis comparing to the sets of  $\mathcal{C}_A, \mathcal{U}_A, Gr_A $ of an ideal, where $\mathbb{N}A$ is not pointed }}}\label{Markovnotpointed}

We have the following theorem, see Figure \ref{relations3}:

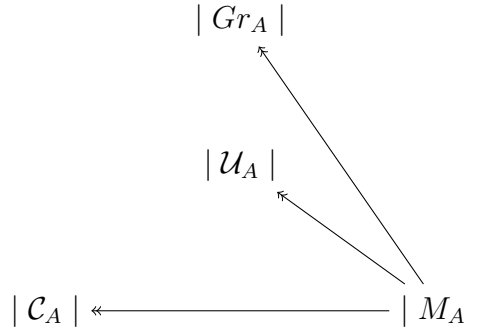
\begin{figure}[h]
\begin{center}
\begin{tikzpicture}[descr/.style={fill=white,inner sep=1.5pt}]
    \matrix (m) [matrix of math nodes,row sep=3em,column sep=3em]
  {
       & \mid Gr_A \mid &  \\
      & \mid \mathcal{U}_A \mid&  \\  
			\mid \mathcal{C}_A \mid&  & \mid M_A	\mid\\		
		};
     \path[-stealth]
    (m-1-2) (m-2-2)
						 (m-3-3) edge[->>] (m-2-2) 
						 (m-1-2) edge[<<-] (m-3-3)
						 (m-3-1) (m-1-2)
						 (m-3-1) edge[<<-] (m-3-3)
						 (m-2-2)  (m-3-1);
\end{tikzpicture}
\caption{Size comparison between all others and a Markov basis}
\label{relations3}
\end{center}
\end{figure}

\begin{thm1}\label{size2} The size of  a Markov basis of $I_A$ cannot
be bounded above by any function on the size of the Graver basis, the universal Gr\"obner basis 
and the set of the circuits of a toric ideal $I_A$.
\end{thm1}
\begin{proof} Let ${\bf a}_1=1$ and ${\bf a}_2=-1$ then $\mathbb{N}A$ is the simplest example of a not pointed semigroup. 

The toric ideal $I_A=\langle x_1x_2-1 \rangle$ is principal and therefore  
$ \mathcal{C}_G = \mathcal{U}_G = Gr_G  =\{x_1x_2-1\}$. 

Let $a, b$ be positive integers and we consider $d=\textrm{gcd}(a,b)$. Then both $x_1^ax_2^a-1$ and $ x_1^bx_2^b-1$ are multiples of $x_1^dx_2^d-1$, 
since $d$ divides both $a, b$. Note also that $d$ can be expressed in the form $ka-lb$ or $lb-ka$ for some non negative integers. Then $x_1^{ka}x_2^{ka}-1-x_1^dx_2^d(x_1^{lb}x_2^{lb}-1)=x_1^dx_2^d-1$ or $x_1^{lb}x_2^{lb}-1-x_1^dx_2^d(x_1^{ka}x_2^{ka}-1)=x_1^dx_2^d-1$
and therefore
$$\langle x_1^ax_2^a-1, x_1^bx_2^b-1\rangle=\langle x_1^dx_2^d-1\rangle.$$ Using induction, the same formula is true for more than two integers. 

Let $q_1,\ldots ,q_s$ be pairwise relative prime integers greater than 1. Let $Q=q_1\cdots q_s$ and $a_i=Q/q_i$. Then 
$\langle x_1^{a_1}x_2^{a_1}-1, x_1^{a_2}x_2^{a_2}-1, \ldots ,x_1^{a_s}x_2^{a_s}-1\rangle= \langle x_1x_2-1\rangle=I_A$ since the greatest common divisor of the
$a_1,\ldots,a_s$ is one and $$\langle x_1^{a_j}x_2^{a_j}-1\ |\  j\not= i\ \text{and}\  1\leq j\leq s\rangle=\langle x_1^{q_i}x_2^{q_i}-1\rangle\not=I_A.$$ Therefore 
$\{x_1^{a_1}x_2^{a_1}-1, x_1^{a_2}x_2^{a_2}-1, \cdots ,x_1^{a_s}x_2^{a_s}-1\}$ is a Markov basis. We have seen for the above example that the size of the Graver basis, the universal Gr\"obner 
basis 
and the set of the circuits of a toric ideal $I_A$ is one while there exist Markov bases of arbitrary large size. 
\end{proof}
This is a classical well known result. It shows  easily an expected behavior of not pointed semigroups with respect to Markov bases. 
For  more examples of not pointed semigroups
and the properties of their Markov bases see \cite{THO2}. 

Combining the results of Theorem \ref{size1} and Theorem \ref{size2} we have the Figure \ref{relations}:

\begin{figure}[h]
\begin{center}
\begin{tikzpicture}[descr/.style={fill=white,inner sep=1.5pt}]
    \matrix (m) [matrix of math nodes,row sep=3em,column sep=3em]
 {
       & \mid {Gr_A} \mid&  \\
      & \mid {\mathcal{U}_A} \mid&  \\  
		\mid {\mathcal{C}_A} \mid&  & \mid{M_A}	\mid\\		
	};
     \path[-stealth]
      (m-3-3) edge[<<->>] (m-1-2)
						 (m-3-3) edge[<<->>] (m-2-2) 
					 (m-3-1) edge[<<->>] (m-3-3);
\end{tikzpicture}
\caption{Comparison of the size of the toric bases}
\label{relations}
\end{center}
\end{figure}
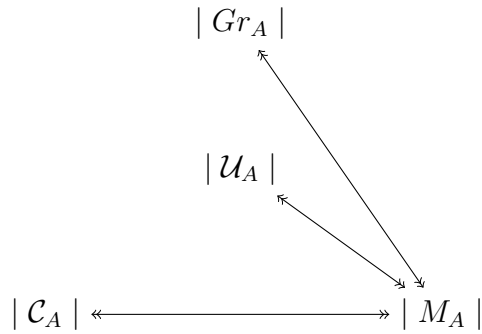

\subsection{On the degree bounds of toric bases}\label{degrees}

\subsubsection{{\bf{Maximal degree of the Graver basis comparing to the sets of  $\mathcal{C}_A, \mathcal{U}_A, M_A $ of an ideal}}}

Let $I_A$ be a toric ideal. For the rest of the paper, we denote by $\deg_{Gr_A},\deg_{\mathcal{U}_A},\deg_{\mathcal{C}_A},\deg_{M_A}$ 
the maximal degree of the elements of the Graver basis, the elements of the universal Gr\"obner basis, 
the circuits and the minimal generators of $I_A$, correspondingly. 

Let $G_1$, $G_2$ be two vertex disjoint graphs,
on the vertices sets $V(G_1)=\{v_1,\ldots,v_s\}$, $V(G_2)=\{u_1,\ldots,u_k\}$
and on the edges sets $E(G_1), E(G_2)$ correspondingly.
We define the {\em sum of the graphs} $G_1,G_2$ on the vertices $v_i,u_j$ as a new graph $G$
formed from their union by identifying the pair of vertices  $v_i,u_j$ to form a single vertex $u$.
The new vertex $u$ is a cut vertex  in the new graph $G$ if both  $G_1$, $G_2$ are not trivial.
We say that we {\em add} to a vertex $v$ of a graph $G_1$ a cycle $S$, to get a graph $G$ if $G$ is the sum of
$G_1$ and $S$ on the vertices $v\in V(G_1)$ and any vertex $u\in S$ correspondingly.

\begin{figure}[h]
\begin{center}
\includegraphics{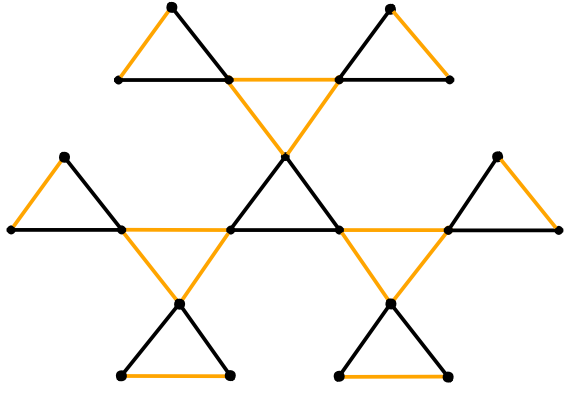}
\caption{The Eulerian trail $w_2^3$ of the graph $G_2^3$}
\label{Figure 44}
\end{center}
\end{figure}

Let $n$ be an odd integer greater than or equal to three. Let $G_0^n$ be a cycle of length $n$.
For $r\ge 0$ we define the graph $G_{r}^n$ inductively on $r$.
$G_{r+1}^n$ is the graph taken from  $G_{r}^n$ by adding to each vertex of
degree two of the graph $G_{r}^n$  a cycle of length $n$. 
Figure \ref{Figure 44} shows the  graph $G_2^3$.

We consider the graphs $G_0^n$ up to $G_{r-1}^n$ as subgraphs of $G_r^n$.
We note that the graph $G_r^n$ is Eulerian since by construction  
it is connected and every vertex has even degree either four if it is also a vertex of $G_{r-1}^n$ or two if it  is not.
Let $w_r^n$ be any closed Eulerian trail of the graph $G_r^n$, i.e.
a trail of the graph which visits every edge of the graph exactly once.

In \cite{TT2} the authors proved that if $w_r^n$ is any closed Eulerian trail of the graph $G_r^n$, then 
the corresponding binomial $B_{w_r^n}$ is an element of the Graver basis of the ideal $I_{G_r^n}$ of degree
 $$\deg(B_{(w_r^n)})=\frac{1}{2}(n+n^2(\frac{(n-1)^r-1}{n-2})). \eqno{(*)}$$
 Moreover, they proved that the maximal degree of the circuits of the above graph is linear on $r$, i.e.
 $$\deg(C_{G_r^n})= n+(2r-1)(n-1), \eqno{(**)}$$ for more see \cite[Proposition 4.3 and Remark 4.4]{TT2}.

 We use these results to prove that the degree of an element in the Graver basis 
$Gr_A$ of a toric ideal $I_{A}$ cannot be bounded above by a polynomial on the maximal 
degree of the elements of the universal Gr\"obner basis $\mathcal{U}_{A}$ and on the maximal degree of the
elements of any Markov basis $M_{A}$ of the ideal too. To prove such a claim, we remark that for the graph $G_r^3$
the corresponding toric ideal is an example of an ideal in which the set of the circuits
and the set of the elements of 
the universal Gr\"obner basis coincide, while the set of the minimal generators of the ideal
is being included strictly on them as we can see in the next proposition. Note that the same result is not true 
if $n\not =3$. To prove the next proposition we define the block tree of a graph $G$. Let $B(G)$ be the {\em block tree} of $G$, the bipartite graph with bipartition
$(\mathbb{B},\mathbb{S})$ where $\mathbb{B}$ is the set of blocks of
$G$ and $\mathbb{S}$ is the set of cut vertices of $G$,
$\{{\mathcal B}, v\}$ is an edge if and only if $v\in {\mathcal B}$. 

\begin{prop1}\label{MU}
We consider the graph $G_r^3$ and let $I_{G_r^3}$ be its corresponding toric ideal. It holds 
$$M_{G_r^3}\subsetneq\mathcal{C}_{G_r^3}=\mathcal{U}_{G_r^3}.$$
\end{prop1}
\begin{proof} Firstly, we will prove that $M_{G_r^3}\subsetneq\mathcal{C}_{G_r^3}$. It is 
enough to find an element of $\mathcal{C}_{G_r^3}$ which is not minimal.
By definition, the graph $G_r^3$ includes the graph $G_1^3$ as a subgraph. By construction, the graph $G_1^3$
 consists of four cycles and let it be
$G_1^3=\{c_1,c_2=(e_1,e_2,e_3),c_3,c_4\}$. Clearly for the walk $w=(c_1,e_1,e_2,c_3)$ 
the corresponding binomial $B_w$ belongs to the $\mathcal{C}_{G_r^3}$. We remark that the 
edge $e_3$ is a chord of $w$ which is a bridge for the graph ${\bf{w}}$, which means that it is not an odd chord of $w$.
It follows from Theorem \ref{minimal} that $B_w$ is not minimal.

Next we will prove that $\mathcal{C}_{G_r^3}=\mathcal{U}_{G_r^3}$. By Proposition \ref{CUG} we know that $\mathcal{C}_{G_r^3}\subseteq\mathcal{U}_{G_r^3}$.
Let $B_w$ be an element of the $\mathcal{U}_{G_r^3}$, where $w$ is a primitive 
even closed walk of $G_r^3$. We will prove that $B_w\in\mathcal{C}_{G_r^3}$. From Theorem \ref{UGB} it follows that the walk $w$ is mixed. 
Therefore every cyclic block of ${\bf w}$ is not pure. But every cycle is a 3-cycle
which means that in every cyclic block of ${\bf w}$ two of the edges are in  ${\bf w}^+$ and one in
${\bf w}^-$ or conversely. That means that each cyclic block has exactly one sink which means also exactly one 
cut vertex. 

Let $n$ be the number of blocks of ${\bf w}$, $c$ the number of cut vertices and $s$ be the number of the cyclic blocks of ${\bf w}$, then $n-s$ is the number
of cut edges. Every cut vertex belongs to exactly two blocks which means that $2c$ is the number of edges in $B({\bf w})$, the block tree 
of ${\bf w}$. The graph $B({\bf w})$ is a tree therefore the number of vertices of $B({\bf w})$ is the number of edges plus one. 
Thus $n+c=2c+1$. It follows that 
there are $n-1$ cut vertices in ${\bf w}$.
Every cut edge of ${\bf w}$ has two cut vertices and every cyclic block we have just proved has only one, therefore
$$s+2(n-s)=2(n-1).$$
Therefore $s=2$ and that means $B_w$ is a circuit since  ${\bf w}$ consists of two odd cycles and possible cut edges, see \cite{Vi}.
\end{proof}

\begin{figure}[h]
\begin{center}
\begin{tikzpicture}[descr/.style={fill=white,inner sep=1.5pt}]
    \matrix (m) [matrix of math nodes,row sep=3em,column sep=3em]
  {
       & \deg_{Gr_A} &  \\
      & \deg_{\mathcal{U}_A} &  \\  
			\deg_{\mathcal{C}_A} &  & \deg_{M_A}	\\		
		};
     \path[-stealth]
    (m-1-2) edge [->>] (m-2-2)
	(m-1-2) edge[->>] (m-3-3)
    (m-3-1) edge[<<-]  (m-1-2);
\end{tikzpicture}
\caption{Comparison of degree bounds I}
\label{relationsA}
\end{center}
\end{figure}
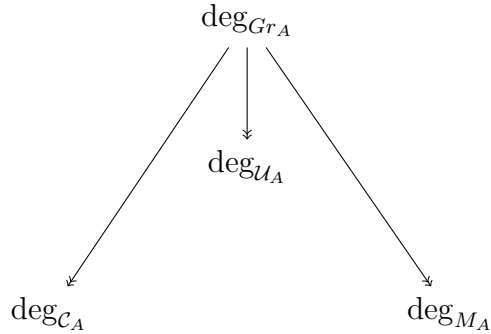

From the above proposition, we have the following theorem, see Figure \ref{relationsA}.

\begin{thm1}\label{TM1} The degrees of the elements in the Graver basis  of a toric ideal $I_A$ cannot
be bounded above by a polynomial on the maximal degree of the circuits, on the maximal degree of 
the elements of a Markov basis 
and on the maximal degree of the elements of the universal Gr\"obner basis of $I_A$.
\end{thm1}
\begin{proof} We consider the graph $G_r^3$. From relation (*) it follows that 
$\deg(B_w)=9\cdot2^{r-1}-3$ for an element of maximal degree of the Graver basis of the ideal. Also, from relation (**) 
we have that $\deg_{\mathcal{C}_{G_r^3}} = 4r+1$. The result follows from Proposition \ref{MU}. 
\end{proof}

\subsubsection{{\bf{Maximal degree of the circuits comparing to the sets of  $Gr_A, \mathcal{U}_A, M_A $ of an ideal}}}

Strongly robust toric ideals are ideals such that the Graver basis is a Markov basis. This implies that there is a unique Markov basis
which is identical with any Gr\"{o}bner basis, thus also with the universal Gr\"{o}bner basis of $I_G$ as well as with the Graver basis.
The only set that we do not have any information is the set of the circuits, except that it is a subset of the Graver basis. The next example shows that there may be huge difference 
between the size of these two sets. 

We recall the definition of the subdivision of a graph.  
 A $k$-subdivision of a graph $G$ is a new graph $S_k(G)$ taken from $G$ by replacing every 
 edge of $G$ by $k$ new edges, where $k\ge 2$. Let $G=(V,E)$ then $$V(S_k(G))=V(G)\sqcup (\bigsqcup \{x(e)_1, \ldots , x(e)_{k-1}\ |\ e\in E \})$$
 $$E(S_k(G))=\bigsqcup \{(u,x(e)_1), (x(e)_1,x(e)_2), \ldots , (x(e)_{k-1},v)\ |\ e=(u,v)\in E \}.$$
The new vertices $x(e)_i$  are all of degree two and this implies that any closed walk of $S_k(G)$ that passes through
any one of $(u,x(e)_1), (x(e)_1,x(e)_2), \cdots , (x(e)_{k-1},v)$ passes through all of them. Therefore there exists
a one to one and onto correspondence between closed walks of $G$ and closed walks of $S_k(G)$. Let $w$ be a closed 
walk of $G$ of length $m$ then the corresponding  closed walk $w_k$ of $S_k(G)$ has length $km$ and has the property 
if $e=(u,v)$ is an edge of the walk $w$ then  $(u,x(e)_1), (x(e)_1,x(e)_2), \cdots , (x(e)_{k-1},v)$ are edges of $w_k$.

\begin{prop1}\label{robust}
Let $G$ be any graph and $S_k(G)$ the $k$ subdivision of $G$. The toric ideal 
$I_{S_k(G)}$ is strongly robust i.e.  
$$M_{S_k(G)}=\mathcal{U}_{S_k(G)}=Gr_{S_k(G)}.$$
\end{prop1}
\begin{proof}
 We will prove that $Gr_{S_k(G)}\subseteq M_{S_k(G)}$. Let 
$B_{w_k}\in Gr_{S_k(G)}$. Then $w_k$ is primitive and also strongly primitive since the distance of any two sinks 
in a cyclic block, if any, is a multiple of $k$. Also in ${\bf w_k}$ there are no chords and certainly no $F_4$.
It follows from 
Theorem \ref{minimal} and Theorem \ref{indispensable} that $B_w$ is minimal and indispensable therefore there exists a unique Markov basis,
thus $B_w\in M_{S_k(G)}$. We note that every Gr\"obner basis contains a Markov basis which implies that 
$Gr_{S_k(G)}\subseteq M_{S_k(G)}\subseteq \mathcal{U}_{S_k(G)}\subseteq Gr_{S_k(G)}.$
\end{proof}

\begin{rem1} {\rm 
The graph $S_k(G)$ is an example
of a graph which is robust, generalized robust (i.e. the union of all Markov bases of the ideal forms its universal Gr\"obner basis) and strongly robust, for more see \cite{BCH,Sull,T}.
We remark that the same result can be proved by using the bouquet algebra of toric ideals developed in \cite{PeThVl}.
Note that the edges $(u,x(e)_1), (x(e)_1,x(e)_2), \cdots , (x(e)_{k-1},v)$ of $S_k(G)$ which correspond to the 
edge $(u,v)$ of $G$ belong to the same bouquet of the toric ideal $I_{S_k(G)}$, since every walk that passes through one of them 
passes through all of them and the corresponding bouquet is mixed since $k>1$. The result follows from \cite[Corollary 4.4]{PeThVl}. 
Note also that if $k$ is odd then a bouquet ideal of $I_{S_k(G)}$ is $I_G$.} 
\end{rem1}

\begin{figure}[h]
\begin{center}
\includegraphics{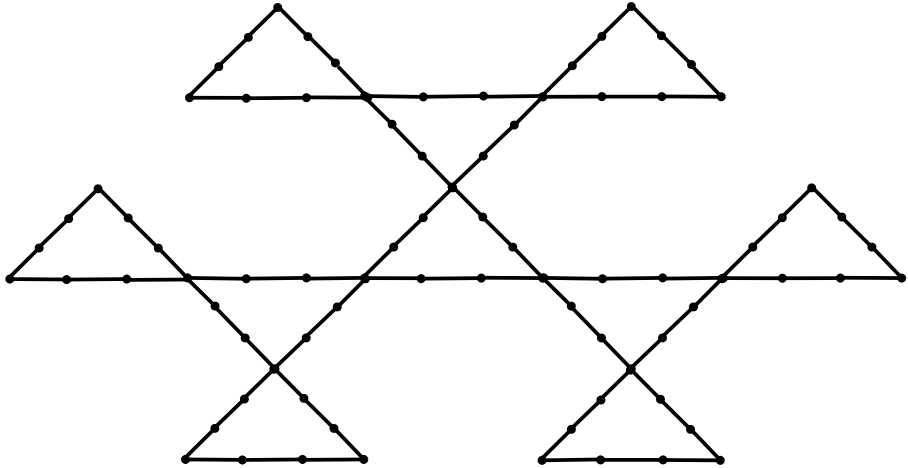}
\caption{A 3-subdivision of the graph $G_2^3$}
\label{Figure!!}
\end{center}
\end{figure}

The graph $S_k(G_{r}^n)$ is defined as a $k$-subdivision
of the graph $G_{r}^n$, where $k$ is odd and $k\geq 3$. For example, we present the graph $S_3(G_{2}^3)$, see Figure \ref{Figure!!}.

\begin{prop1}\label{UM}
We consider the graph $S_k(G_{r}^n)$ and let $I_{S_k(G_{r}^n)}$ be its corresponding toric ideal. 
For any closed Eulerian trail $(w_k)_r^n$ of the graph, we have that $B_{(w_k)_r^n}\in Gr_{S_k(G_r^n)}$. 
\end{prop1}
\begin{proof}
Let $(w_k)_r^n$ be a closed Eulerian trail of the graph $S_k(G_{r}^n)$, then  
$w_r^n$ is a closed Eulerian trail of the graph $G_{r}^n$. 
The graph of the walk $w_r^n$ is  $G_{r}^n$ and the graph of the walk
$(w_k)_r^n$ is $S_k(G_{r}^n)$. In \cite[Proposition 4.1.]{TT2}, the authors
proved for the graph $G_{r}^n$ that $B_{w_r^n}$ is an element of the Graver basis of $I_{G_r^n}$.
The graph $S_k(G_{r}^n)$  is not biconnected and every 
block is a cycle. Also, if $v$ is a cut vertex of ${S_k(G_r^n)}$, it is also a cut vertex of ${G_r^n}$. The vertex $v$
separates the graph ${G_r^n}$ in two parts and let $m_1,m_2$ be the total number of the 
edges of the blocks that are cycles in each part correspondingly. Since the binomial $B_{w_r^n}$ 
is an element of the Graver basis, we remark that $m_1,m_2$ are odd. Obviously, the two parts of the graph ${S_k(G_r^n)}$,
which separates the vertex $v$ have $km_1,km_2$ total number of edges correspondingly, i.e. an odd number. 
By Theorem \ref{primitive-graph}, it follows that
the binomial $B_{(w_k)_r^n}$ belongs to the Graver basis of $I_{S_k(G_{r}^n)}$. 

\end{proof}

We consider the graph $B(S_k(G_{r}^n))$ to be the block tree of $S_k(G_{r}^n)$.
Let ${\mathcal B}_k, {\mathcal B}_i, {\mathcal B}_l$ be blocks of a graph $S_k(G_{r}^n)$. We call the block 
${\mathcal B}_i$ an \emph{internal block}
 of ${\mathcal B}_k,{\mathcal B}_l$, if ${\mathcal B}_i$ is an internal vertex in the unique path defined 
 by ${\mathcal B}_k,{\mathcal B}_l$ in the block tree $B(S_k(G_{r}^n))$.
Every path of the graph $S_k(G_{r}^n)$ from the block ${\mathcal B}_k$ to the block ${\mathcal B}_l$ passes
through every internal  block of ${\mathcal B}_k,{\mathcal B}_l$.
The path has vertices at least the cut vertices of $S_k(G_{r}^n)$ which are vertices in
 the  path $({\mathcal B}_k,\ldots,{\mathcal B}_l)$
  in $B(S_k(G_{r}^n))$.
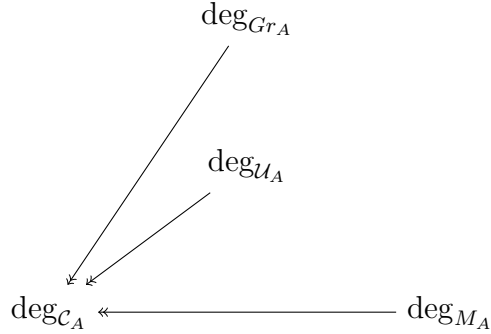
\begin{figure}[h]
\begin{center}
\begin{tikzpicture}[descr/.style={fill=white,inner sep=1.5pt}]
    \matrix (m) [matrix of math nodes,row sep=3em,column sep=3em]
  {
       & \deg_{Gr_A} &  \\
      & \deg_{\mathcal{U}_A} &  \\  
			\deg_{\mathcal{C}_A} &  & \deg_{M_A}	\\		
		};
     \path[-stealth]
    (m-1-2) edge [->>] (m-3-1)
	(m-2-2) edge[->>] (m-3-1)
    (m-3-1) edge[<<-]  (m-3-3);
\end{tikzpicture}
\caption{Comparison of degree bounds II}
\label{relationsB}
\end{center}
\end{figure}

\begin{thm1}\label{TM2} The degrees of the elements of a Markov basis, the degrees of the elements
of the universal Gr\"obner basis and the degrees of the elements of the Graver basis of a toric ideal $I_A$, 
cannot be bounded above by a polynomial on the degrees of the circuits of the ideal $I_A$.
\end{thm1}
\begin{proof}
We consider the graph 
$S_k(G_{r}^n)$ and let $I_{S_k(G_{r}^n)}$ be its corresponding toric ideal. 
Let $(w_k)_r^n$ be an
Eulerian trail of the graph $S_k(G_{r}^n)$. By Proposition \ref{UM} we have that
$B_{(w_k)_r^n}$ is an element of the Graver basis of $I_{S_k(G_{r}^n)}$. 
For the Eulerian trail $w_r^n$ of the graph $G_{r}^n$
from the relation (*) we have that $\deg(B_{w_r^n})=\frac{1}{2}(n+n^2(\frac{(n-1)^r-1}{n-2})). $ 
Every edge of the walk $w_r^n$ corresponds to $k$ edges of the walk $(w_k)_r^n$ therefore 
 $$\deg(B_{(w_k)_r^n})=\frac{k}{2}(n+n^2(\frac{(n-1)^r-1}{n-2})).$$
Which means that there exists an element in the Graver basis of $I_{S_k(G_{r}^n)}$ 
whose degree  is exponential in $r$. By Proposition \ref{robust}, it happens as well for the elements 
of the Markov basis and of the universal Gr\"obner basis of the ideal.\par

Let $B_w$ be a circuit of $I_{S_k(G_{r}^n)}$. The graph $S_k(G_{r}^n)$ has no even 
cycles and therefore the subgraph corresponding to a circuit  
consists of two different  odd cycles joined by a path, see Theorem \ref{circuit}.  We remark that 
every cycle of the graph $S_k(G_{r}^n)$  has length $kn$ and it is a block. We claim that a path between two 
blocks ${\mathcal B}_1,{\mathcal B}_2$ of $S_k(G_{r}^n)$ has length at most $(2r-1)k(n-1)$. 
Each such path passes through all internal blocks of ${\mathcal B}_1,{\mathcal B}_2$ and no other 
and has at most $k(n-1)$ common edges  with 
every one of them. We denote by $d({\mathcal B}_1,{\mathcal B}_2)$  the number of the internal blocks of
${\mathcal B}_1,{\mathcal B}_2$. From \cite[Lemma 4.2.]{TT2} we know that $d({\mathcal B}_1,{\mathcal B}_2)\leq 2r-1.$
Therefore the path has at most length $d({\mathcal B}_1,{\mathcal B}_2)\cdot k(n-1)\leq (2r-1)k(n-1).$ 
Thus the corresponding circuit has
degree at most $ kn+(2r-1)k(n-1)$ which is linear on $r$.\par Therefore the  degree of an element in the Graver basis, 
of the Markov basis and of the universal Gr\"obner basis
of the toric ideal $I_{A_{S_k(G_{r}^n)}}$ cannot
be bounded above by a polynomial on the maximal degree of a circuit.
\end{proof}
\begin{rem1}{\rm Note that one can use the main result of \cite{TT2} and Lawrence liftings to prove the same result. Lawrence liftings provide a technique 
of constructing toric ideals which are strongly robust. For more on strongly robust ideals and a generalization of the Lawrence lifting technique one can
read \cite{PeThVl}. 
}
\end{rem1}

\subsubsection{{\bf{Maximal degree of the Markov basis comparing to the sets of  $\mathcal{C}_A, \mathcal{U}_A, Gr_A $ of an ideal}}}

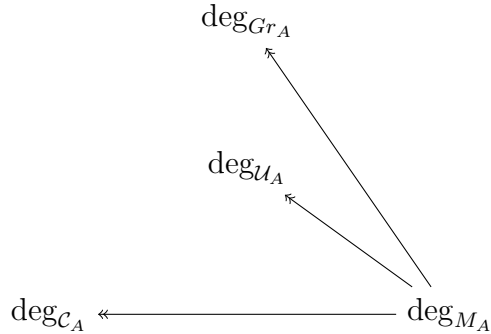
\begin{figure}[h]
\begin{center}
\begin{tikzpicture}[descr/.style={fill=white,inner sep=1.5pt}]
    \matrix (m) [matrix of math nodes,row sep=3em,column sep=3em]
  {
       & \deg_{Gr_A} &  \\
      & \deg_{\mathcal{U}_A} &  \\  
			\deg_{\mathcal{C}_A} &  & \deg_{M_A}	\\		
		};
     \path[-stealth]
    (m-3-3) edge [->>] (m-2-2)
	(m-3-1) edge[<<-] (m-3-3)
    (m-3-3) edge[->>]  (m-1-2);
\end{tikzpicture}
\caption{Comparison of degree bounds III}
\label{relationsD}
\end{center}
\end{figure}

We recall the example of Subsection \ref{Markovnotpointed}, where we are in the case that $\mathbb{N}A$ is not pointed.
By looking at the maximal degree of the elements of the different bases we remark that the unique element in
the Graver basis, the universal Gr\"obner 
basis 
and the set of the circuits of the toric ideal $I_A$ has degree 2 while the maximal degree of an element in a Markov basis can be arbitrary high.
Therefore we have also the following theorem, see Figure \ref{relationsD}:
\begin{thm1}\label{TM4} The degrees of the elements of a Markov basis of $I_A$ cannot
be bounded above by any function on the maximal degree of the elements of the Graver basis, on the maximal degree of the elements of the universal Gr\"obner basis 
and on the maximal degree of the of the circuits of a toric ideal $I_A$.
\end{thm1}

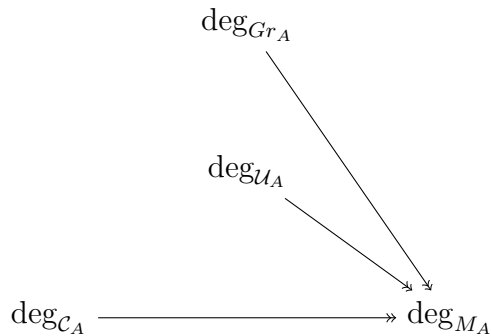
\begin{figure}[h]
\begin{center}
\begin{tikzpicture}[descr/.style={fill=white,inner sep=1.5pt}]
    \matrix (m) [matrix of math nodes,row sep=3em,column sep=3em]
  {
       & \deg_{Gr_A} &  \\
      & \deg_{\mathcal{U}_A} &  \\  
			\deg_{\mathcal{C}_A} &  & \deg_{M_A}	\\		
		};
     \path[-stealth]
    (m-1-2) edge [->>] (m-3-3)
	(m-2-2) edge[->>] (m-3-3)
    (m-3-3) edge[<<-]  (m-3-1);
\end{tikzpicture}
\caption{Comparison of degree bounds IV}
\label{relationsC}
\end{center}
\end{figure}

Returning in the case that $\mathbb{N}A$ is pointed, we consider the complete graph $G$. We note that a graph is complete if each pair of two distinct vertices of $G$, is connected by an edge.
The complete graph with $n$ vertices is denoted by $K_n$. On this example we have the following theorem, see Figure \ref{relationsC}:

\begin{thm1}\label{TM3} The degrees of the elements of the Graver basis, the degrees of the elements of the universal Gr\"obner basis and the degrees of the circuits of the ideal cannot be bounded above by any function on the maximal degrees
of the elements of a Markov basis of the ideal.
\end{thm1}
\begin{proof} In \cite[Proposition 4.1]{TT1}, the authors proved that 
the largest degree $d_n$ of a binomial in the Graver basis (and in the universal Gr\"{o}bner basis) 
for $I_{K_n}$ is $d_n=n-2$, for $n\geq 4$ and it is attained by a circuit, see also \cite{LST}. Furthermore, 
in \cite{LST}, the authors proved that the maximal degree of the elements of a Markov basis of the ideal $I_{K_n}$ is two.
The result follows.
\end{proof}

Combining the results of Theorem \ref{TM1}, Theorem \ref{TM2}, Theorem \ref{TM4} and Theorem \ref{TM3} we conclude them in Figure \ref{Finalrelations}:
 
 \begin{figure}[h]
\begin{center}
\begin{tikzpicture}[descr/.style={fill=white,inner sep=1.5pt}]
    \matrix (m) [matrix of math nodes,row sep=3em,column sep=3em]
  {
       & \deg_{Gr_A} &  \\
      & \deg_{\mathcal{U}_A} &  \\  
			\deg_{\mathcal{C}_A} &  & \deg_{M_A}	\\		
		};
     \path[-stealth]
    (m-3-3) edge [->>] (m-2-2)
     (m-2-2) edge [->>] (m-3-3)
      (m-2-2) edge [->>] (m-3-1)
	(m-3-1) edge[<<-] (m-3-3)
		(m-3-3) edge[<<-] (m-3-1)
		 (m-1-2) edge[->>]  (m-2-2)
		  (m-1-2) edge[->>]  (m-3-1)
		    (m-1-2) edge[->>]  (m-3-3)
    (m-3-3) edge[->>]  (m-1-2);
\end{tikzpicture}
\caption{Comparison of degree bounds of toric bases}
\label{Finalrelations}
\end{center}
\end{figure}
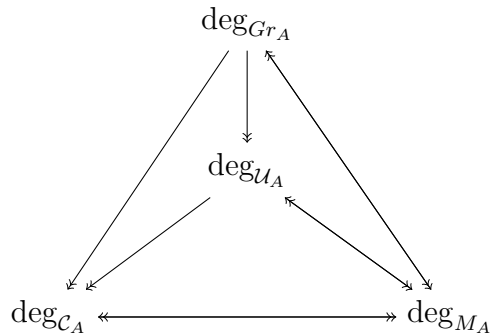

\section{Conclusion}

The main results of this manuscript are presented in Figure \ref{relations} and Figure \ref{Finalrelations}. For the maximal degree of the elements of the sets: Graver basis, universal Gr\"{o}bner basis, a Markov basis and  the set of the circuits of a toric ideal
we managed to prove that  for any $A,B$ of these bases such that $A\not \supseteq B$ in general, there is no polynomial
on the maximal degree of the elements of $A$  which bounds the
maximal degree of the elements of $B$ correspondingly, for any possible
combination. For the size we were able to provide the corresponding
theorems in all cases except three, see Figure \ref{relations}. Since maximal degree of a set can be considered as a measure of the size of this set
we thought that the counterexamples used to prove the theorems in Subsection \ref{degrees}, 
may be also counterexamples in proving the missing theorems about the bounds of the size. But the form of the elements in these bases  
were so complicated that we were not be able to provide formulas for the actual size of these sets. Theorems \ref{circuit},\ref{primitive-graph},\ref{minimal},\ref{UGB} 
suggest to use graphs like $G^n_r$ to produce counterexamples. The advantage of graphs like  $G^n_r$ is that they 
increase the number of blocks exponentially on $r$ but make
computation of all the elements even in the Graver basis very complicated.
It seems to us that to find such counterexamples one has to move to a different class of toric ideals than the toric ideals of graphs.

Although most of the examples used in this article are toric ideals of graphs, one can use the theory of stable toric ideals 
and the generalized Lawrence matrices, developed in \cite{PeThVl}, to produce examples of more general toric ideals
that they have exactly the same properties, since stable toric ideals preserve the size of toric bases and choosing all the vectors $c_B$
to have the same 1-norm $k$ all the degrees are multiplied with the constant $k$, for more details and examples see \cite{PeThVl}.

In \cite{RYCH} it was proved that for any toric ideal $I_G$ of a graph $G$ the degree of any element of the Graver basis of $I_G$
is bounded above by an exponential function of the maximal degree of a circuit. It is an interesting problem if this is true
for any toric ideal. To prove something like that one needs a better understanding of circuits, Markov basis, universal Gr\"{o}bner
and Graver basis for general toric ideals or equivalently, in the case of pointed affine semigroups, for toric ideals of hypergraphs,
see \cite{PeThVl2}.

\end{document}